# First Passage through a Continuous Barrier: Pathwise Decomposition, Random-Time Structure, and Compensators

Tristan Guillaume [1]

[1] CY Cergy Paris Université, Laboratoire Thema, Cergy, France
Correspondence: Tristan Guillaume, CY Cergy Paris Université, Laboratoire Thema,
33 boulevard du port, F-95011 Cergy-Pontoise Cedex, France.
Tel : 33-6-12-22-45-88.
Email : tristan.guillaume@cyu.fr

## Abstract

Let $\tau$ be the first-passage time of a continuous barrier by a càdlàg adapted process. We show that $\tau$ admits a canonical fourfold pathwise decomposition into continuous contact, contact from the left followed by an upward jump, exact hit by jump, and strict overshoot by jump from below. This refinement is more informative than the classical contact-versus-overshoot dichotomy for random-time purposes, because it separates modes with different predictability properties. In particular, the left-contact component always defines an accessible stopping time and becomes predictable under a no-premature-left-contact condition, which we prove to be both sufficient and necessary for the canonical running-supremum announcing sequence to work. On the gap side, under a structural exclusion of predictable gap-crossings, the corresponding restricted time is totally inaccessible. In the semimartingale setting, we obtain a sharp compensator criterion for the predictable-side condition, explicit compensator formulas for the jump-driven crossing modes, and a decomposition of the compensator of the default indicator into its predictable jump part and continuous part. As an application, for a mean-reverting affine jump-diffusion with upward exponential jumps, we derive the boundary-value problem governing the overshoot mode, prove that the differentiated third-order ODE is equivalent to the original problem only when a boundary compatibility condition is retained, and establish verification and uniqueness for the discounted problem. This yields an explicit Green-Volterra representation, a first-order small-$q$ expansion expansion, and, in the undiscounted case, closed formulas for the overshoot and creeping probabilities.

# 1    Introduction

Barrier-crossing times are among the most basic examples of random times. They arise in fluctuation theory, stochastic calculus, mathematical finance, and credit-risk modelling. For càdlàg models, however, the phrase "first passage" conceals several genuinely different geometric mechanisms: the process may arrive continuously, reach the barrier from the left and then jump, hit the barrier exactly by a jump, or overshoot it strictly. The purpose of this paper is to isolate that geometry and relate it to the nature of the crossing time itself.

For Lévy processes and related Markov models, first passage is usually discussed in terms of contact (creeping) versus passage by jump overshoot [9,18,19]. In the general theory of stochastic processes, by contrast, the relevant taxonomy distinguishes predictable, accessible, and totally inaccessible random times [12,21]. Our starting point is that these two viewpoints are related but not identical, and that a finer pathwise decomposition is needed before the random-time structure becomes visible.

Let $X$ be a càdlàg adapted real-valued process on a filtered probability space satisfying the usual conditions, and let $b$ be a continuous real-valued function. We consider the first-passage time

$$\tau := \inf\{t \geq 0 : X_t \geq b(t)\}. \tag{1}$$

Writing $Y_t := X_t - b(t)$, the continuity of $b$ and the càdlàg property of $X$ imply that, on $\{\tau < \infty\}$,

$$Y_{\tau-} \leq 0 \leq Y_\tau. \tag{2}$$

This yields a canonical pathwise classification of first passage into four mutually exclusive modes:

$$\{Y_{\tau-} = 0, Y_\tau = 0\}, \quad \{Y_{\tau-} = 0, Y_\tau > 0\}, \quad \{Y_{\tau-} < 0, Y_\tau = 0\}, \quad \{Y_{\tau-} < 0, Y_\tau > 0\}. \tag{3}$$

These correspond respectively to continuous contact, contact from the left followed by an upward jump, exact hit by jump, and strict overshoot by jump from below.

This decomposition refines the familiar contact/overshoot dichotomy in exactly the direction needed for random-time questions. In particular, the event $\{Y_{\tau-} = 0\}$ splits into a genuinely continuous-contact piece and an exact-hit-by-jump piece, while the event $\{Y_{\tau-} < 0\}$ records whether the barrier was already visible from the left just before passage. The decisive structural distinction is therefore not merely whether the terminal value lands on the barrier, but whether first passage occurs from left contact or from a strict left gap.

That distinction governs the first group of results. On the left-contact event, the restricted left-contact time is always accessible. Under a no-premature-left-contact condition, the canonical sequence built from the running supremum of $Y$ announces this time, and this condition is sharp: it is exactly what makes that canonical announcing mechanism work. On the gap side, first passage can only occur by jump; under a natural exclusion of predictable gap-crossings, the corresponding restricted gap time is totally inaccessible.

In the semimartingale setting we go further in two directions. First, we characterize the no-premature-left-contact condition by a compensator criterion and derive exact compensator formulas for the jump-driven crossing modes. Second, Theorem 2.19 proves that the compensator of the default indicator splits canonically into a predictable jump contribution generated by left contact and a continuous contribution generated by gap crossing. This gives a direct bridge from barrier-crossing geometry to the hazard-process viewpoint used in credit-risk modelling.

On the analytic side, we then specialise to a mean-reverting affine jump-diffusion with upward exponential jumps. For the overshoot mode we derive the governing integro-differential equation, show that the differentiated third-order ODE is equivalent to the original problem only after a boundary compatibility condition is imposed, and prove verification and uniqueness for the discounted boundary-value problem in the class of bounded classical solutions. This yields an explicit Green–Volterra representation in parabolic-cylinder functions together with a first-order small-$\theta$ expansion whose leading correction is the overshoot-time moment. In the undiscounted limit, the representation collapses to a closed formula for the overshoot probability, from which the creeping probability follows by complementarity.

These results should be positioned against several strands of literature. In fluctuation theory for Lévy processes, the distinction between creeping and passage by jump is classical [9,18,19], and the recent work of Chaumont and Pellas [7] shows how delicate the creeping question already becomes for general curves. On the analytic side, explicit first-passage transforms are available for a number of special jump models, including the double-exponential jump-diffusion of Kou and Wang [16], piecewise exponential Markov processes with two-sided jumps [13], Ornstein–Uhlenbeck processes with jumps [22], and hyper-exponential or mixed-exponential jump models [5,8,25].

The present contribution is different in three respects. First, the structural framework of Sections 2–3 is pathwise and filtration-theoretic: it isolates the accessible/predictable side and the totally inaccessible side of first passage, together with a decomposition of the default indicator compensator, rather than only transforms of the total passage time. Second, the affine analysis is not a routine modification of Kou–Wang [16]: their model is spatially homogeneous and leads to constant-coefficient ODEs, so neither the parabolic-cylinder reduction, nor the boundary compatibility issue behind Lemma 4.4, nor the Green–Volterra representation appears there. Third, Novikov [22] studies transforms of the total first-passage time for Ornstein–Uhlenbeck processes with jumps, but does not separate creeping and overshoot into distinct probabilistic objects with their own boundary-value problems.

From the viewpoint of credit risk, first-passage structural models go back at least to Black and Cox [4], and jump extensions such as Zhou [26] show how sudden default can arise even in a structural framework. On the reduced-form and incomplete-information side, the hazard-process interpretation of default is developed in [1,6,10,11,14]. What seems to be missing there is an endogenous decomposition of the default compensator directly from the barrier-crossing geometry of the underlying càdlàg process. Theorem 2.19 supplies precisely such a decomposition in the full-information first-passage setting.

The paper is organised as follows. Section 2 introduces the setting, notation, and pathwise decomposition; develops the accessible, predictable, and totally inaccessible components; and concludes with a compensator decomposition of the default indicator. Section 3 shows how the fourfold decomposition collapses to the contact-versus-overshoot dichotomy, adds the semimartingale criteria, and gives illustrative examples. Section 4 introduces the affine jump-diffusion model and derives the general overshoot equation. Section 5 treats the discounted mean-reverting case, including the verification and uniqueness theorem, the Green–Volterra representation, and the first-order small-θ expansion. Section 6 gives the undiscounted overshoot probability and the creeping formula. Section 7 concludes.

## 2 Pathwise decomposition and random-time consequences

### 2.1 Setting and notation

Let $(\Omega, \mathcal{F}, \mathbb{F}, \mathbb{P})$ satisfy the usual conditions, where $\mathbb{F} = (\mathcal{F}_t)_{t \geq 0}$. Let $X = (X_t)_{t \geq 0}$ be a càdlàg adapted real-valued process, and let $b: [0, \infty) \to \mathbb{R}$ be continuous. Define

$$Y_t := X_t - b(t), \qquad t \geq 0, \tag{4}$$

and

$$\tau := \inf\{t \geq 0 : Y_t \geq 0\} = \inf\{t \geq 0 : X_t \geq b(t)\}, \tag{5}$$

with the convention $\inf \varnothing = \infty$. Since $X$ is càdlàg adapted and $b$ is continuous deterministic, $Y$ is again càdlàg and adapted. By the début theorem [3], $\tau$ is a stopping time.

We write

$$L := \{\tau < \infty, \ Y_{\tau-} = 0\}, \qquad G := \{\tau < \infty, \ Y_{\tau-} < 0\}. \tag{6}$$

Thus $L$ is the event of *contact from the left*, while $G$ is the event of *gap from the left*. We refine $L$ and $G$ by examining the terminal value $Y_\tau$:

$$C_0 := \{\tau < \infty, \ Y_{\tau-} = 0, \ Y_\tau = 0\}, \qquad \text{(continuous contact)}, \tag{7}$$

$$C_+ := \{\tau < \infty, \ Y_{\tau-} = 0, \ Y_\tau > 0\}, \qquad \text{(left contact, then jump up)}, \tag{8}$$

$$J_0 := \{\tau < \infty, \ Y_{\tau-} < 0, \ Y_\tau = 0\}, \qquad \text{(exact hit by jump)}, \tag{9}$$

$$J_+ := \{\tau < \infty,\ Y_{\tau-} < 0,\ Y_\tau > 0\}, \qquad \text{(strict overshoot by jump)}. \tag{10}$$

By construction, $L = C_0 \cup C_+$ and $G = J_0 \cup J_+$, with each union disjoint. Finally, define the restricted random time associated with left contact by

$$\tau_L(\omega) := \begin{cases} \tau(\omega), & \omega \in L, \\ \infty, & \omega \notin L. \end{cases} \tag{11}$$

## 2.2 Pathwise decomposition at first passage

**Proposition 2.1** (Basic pathwise inequalities). *On the event $\{\tau < \infty\}$, one has*

$$Y_{\tau-} \leq 0 \leq Y_\tau, \qquad \text{equivalently:} \qquad X_{\tau-} \leq b(\tau) \leq X_\tau. \tag{12}$$

*Proof.* For every $s < \tau$, the definition of $\tau$ gives $Y_s < 0$. Passing to the limit along $s \uparrow \tau$ yields $Y_{\tau-} \leq 0$. Since $\tau$ is the first time at which $Y$ enters $[0, \infty)$, one has $Y_\tau \geq 0$. The second inequality follows by adding $b(\tau)$. □

**Proposition 2.2** (Fourfold decomposition). *On $\{\tau < \infty\}$, exactly one of the events $C_0, C_+, J_0, J_+$ occurs. Equivalently,*

$$\{\tau < \infty\} = C_0 \,\dot\cup\, C_+ \,\dot\cup\, J_0 \,\dot\cup\, J_+, \tag{13}$$

*where $\dot\cup$ denotes disjoint union. Moreover, $L = C_0 \,\dot\cup\, C_+$ and $G = J_0 \,\dot\cup\, J_+$.*

*Proof.* By Proposition 2.1, on $\{\tau < \infty\}$ one has $Y_{\tau-} \leq 0 \leq Y_\tau$. Hence either $Y_{\tau-} = 0$ or $Y_{\tau-} < 0$, and independently either $Y_\tau = 0$ or $Y_\tau > 0$. This yields exactly the four alternatives, and their mutual exclusivity is immediate. □

*Remark 2.3* The event $\{Y_\tau = 0\}$ is the union $C_0 \cup J_0$. Thus "contact at time $\tau$" and "continuous contact" are not the same notion. Likewise, $\{Y_\tau > 0\} = C_+ \cup J_+$, so "overshoot at time $\tau$" and "strict gap from the left" are different notions. For random-time purposes, the more basic dichotomy is $\{\tau < \infty\} = L \,\dot\cup\, G$, that is, contact from the left versus gap from the left.

*Remark 2.4* (Information available just before $\tau$). Because $Y_{\tau-}$ is predictable, $Y_{\tau-}$ is $\mathcal{F}_{\tau-}$–measurable. Since $\{\tau < \infty\} = \bigcup_{n \geq 1}\{\tau < n\} \in \mathcal{F}_{\tau-}$, it follows that $L = \{\tau < \infty, Y_{\tau-} = 0\}$ and $G = \{\tau < \infty, Y_{\tau-} < 0\}$ belong to $\mathcal{F}_{\tau-}$.

By contrast, the finer refinements $C_0 = L \cap \{Y_\tau = 0\}$, $C_+ = L \cap \{Y_\tau > 0\}$, $J_0 = G \cap \{Y_\tau = 0\}$, and $J_+ = G \cap \{Y_\tau > 0\}$ depend on the terminal value $Y_\tau$, which is $\mathcal{F}_\tau$–measurable but need not be $\mathcal{F}_{\tau-}$–measurable in general when jumps at $\tau$ are possible. Thus the primary $L$-versus-$G$ dichotomy is the one visible just before $\tau$, whereas the secondary split by the terminal value is typically revealed only at $\tau$.

## 2.3 The predictable side

Define the running supremum process

$$S_t := \sup_{0 \leq s \leq t} Y_s, \qquad t \geq 0, \tag{14}$$

and, for each $n \geq 1$,

$$\sigma_n := \inf\{t \geq 0 : S_t \geq -1/n\}. \tag{15}$$

The unconditional claim that $\tau_L$ is always predictable is false in general: if the path makes an earlier left contact with the barrier and then jumps back down before the actual first-passage time, then the times $\sigma_n$ need not converge to $\tau$. The following condition excludes precisely that pathology.

**Assumption 2.5** (No premature left contacts). For every finite $t \geq 0$,

$$\sup_{0 \leq s \leq t} Y_s < 0 \qquad \text{on } \{\tau > t\}. \tag{16}$$

**Proposition 2.6** *Assumption 2.5 holds if, pathwise, $Y_{u-} < 0$ for every finite $u < \tau$. In particular, Assumption 2.5 holds whenever $Y$ is continuous.*

*Proof.* Fix a sample point $\omega$, and suppose that $Y_{u-}(\omega) < 0$ for every finite $u < \tau(\omega)$. We show that Assumption 2.5 holds along this path. Let $t < \tau(\omega)$. Since $t$ lies strictly before the first-passage time, one has $Y_s(\omega) < 0$ for all $0 \leq s \leq t$. Thus the only possible obstruction to $\sup_{0 \leq s \leq t} Y_s(\omega) < 0$ would be that the supremum is equal to 0. Assume, for contradiction, that $\sup_{0 \leq s \leq t} Y_s(\omega) = 0$. Then there exists a sequence $(s_n) \subset [0, t]$ such that $Y_{s_n}(\omega) \to 0$. Since $[0, t]$ is compact, some subsequence of $(s_n)$ converges to a point $u \in [0, t]$. Passing to a further subsequence if necessary, we may suppose that the sequence is monotone. There are then two possibilities. If $s_n \uparrow u$, then by definition of the left limit, $Y_{s_n}(\omega) \to Y_{u-}(\omega)$. Since also $Y_{s_n}(\omega) \to 0$, it follows that $Y_{u-}(\omega) = 0$. But $u \leq t < \tau(\omega)$, so $u < \tau(\omega)$, contradicting the hypothesis that $Y_{u-}(\omega) < 0$ for every finite $u < \tau(\omega)$. If instead $s_n \downarrow u$, then by right-continuity of $Y$, $Y_{s_n}(\omega) \to Y_u(\omega)$. Again $Y_{s_n}(\omega) \to 0$, hence $Y_u(\omega) = 0$. But $u \leq t < \tau(\omega)$, so $u < \tau(\omega)$, and by definition of $\tau$ this is impossible, since $Y_s(\omega) < 0$ for all $s < \tau(\omega)$. Both cases lead to a contradiction. Therefore, $\sup_{0 \leq s \leq t} Y_s(\omega) < 0$. Since $t < \tau(\omega)$ was arbitrary, Assumption 2.5 holds along the path. For the final claim, if $Y$ is continuous, then $Y_{u-} = Y_u < 0$ for every $u < \tau$, so the first part applies immediately. □

**Proposition 2.7**. *The process $S$ is adapted, càdlàg, and nondecreasing. Consequently, each $\sigma_n$ and each $\rho_n := \sigma_n \wedge n$ is a stopping time.*

*Proof.* Fix $t \geq 0$. Since $Y$ is adapted, for each $q \in \mathbb{Q} \cap [0, t]$ the random variable $Y_q$ is $\mathcal{F}_t$-measurable. Because $Y$ is càdlàg, the supremum over $[0, t]$ agrees pathwise with the supremum over a dense countable set of rational times $S_t = sup_{q \in \mathbb{Q} \cap [0,t]} Y_q$. The right-hand side of the previous inequality is the supremum of countably many $\mathcal{F}_t$-measurable random variables, and is therefore $\mathcal{F}_t$-measurable. Hence $S$ is adapted.

The non-decreasing nature of $S$ is clear. We now verify that $S$ is càdlàg. Since $S$ is nondecreasing, left limits exist automatically. It remains to check right-continuity. Let $u \downarrow t$. Then, $S_u = \max(S_t, sup_{t<s\leq u}Y_s)$. Because $Y$ is right-continuous at $t$, $sup_{t<s\leq u}Y_s \to Y_t$ as $u \downarrow t$. Since $Y_t \leq S_t$, it follows that $S_u \to S_t$. Thus $S$ is right-continuous. Therefore $S$ is càdlàg.

Now, since $\sigma_n$ is the first hitting time of the closed set $[-1/n, \infty)$ by a càdlàg adapted process, $\sigma_n$ is therefore a stopping time. Finally, since the minimum of a stopping time and a deterministic time is again a stopping time, $\rho_n$ is a stopping time. □

**Proposition 2.8** *The restricted random time $\tau_L$ is a stopping time.*

*Proof.* Fix $t \geq 0$. By definition of $\tau_L$, the event $\{\tau_L \leq t\}$ can occur only on $L$, because on $L^c$ one has $\tau_L = \infty$. On the event $L$, however, $\tau_L = \tau$. Therefore, $\{\tau_L \leq t\} = L \cap \{\tau \leq t\}$. Now $L \in \mathcal{F}_\tau$. By the defining property of $\mathcal{F}_\tau$-measurability, this implies that for every $t \geq 0$,

$$L \cap \{\tau \leq t\} \in \mathcal{F}_t.$$

Hence, $\{\tau_L \leq t\} \in \mathcal{F}_t$ for every $t \geq 0$. □

**Proposition 2.9** (Accessibility of the left-contact component). *The restricted random time $\tau_L$ is an accessible stopping time.*

*Proof.* On $\{\tau_L < \infty\}$ one has $\tau_L = \tau$ and, by definition of $L$, $Y_{\tau_L -} = Y_{\tau -} = 0$. Hence the graph of $\tau_L$ is contained in the predictable set $\{(\omega, t) : Y_{t-}(\omega) = 0\}$. A standard theorem of the general theory states that every predictable set is a countable union of graphs of predictable stopping times. Therefore the graph of $\tau_L$ is contained in a countable union of graphs of predictable stopping times, which is exactly the definition of accessibility. □

**Theorem 2.10** (Predictability of the left-contact component). *Assume 2.4. Define $\rho_n := \sigma_n \wedge n$ for $n \geq 1$. Then $(\rho_n)_{n\geq 1}$ is an increasing sequence of stopping times such that*

$$\rho_n < \tau_L \text{ on } \{\tau_L < \infty\}, \qquad \rho_n \uparrow \tau_L \text{ pointwise.} \tag{17}$$

*In particular, $\tau_L$ is predictable.*

*Proof.* Each $\rho_n$ is a stopping time by Proposition 2.7. Since $S$ is nondecreasing and $-1/(n+1) > -1/n$, one has $\sigma_n \leq \sigma_{n+1}$, so $\rho_n \leq \rho_{n+1}$.

If $\omega \in L$, then $\tau(\omega) < \infty$, $Y_{\tau-}(\omega) = 0$, and $\tau_L(\omega) = \tau(\omega)$. By Assumption 2.5, for every $t < \tau(\omega)$, $S_t(\omega) < 0$. On the other hand, since $Y_{\tau-}(\omega) = 0$, there exists a sequence $s_k \uparrow \tau(\omega)$ such that $Y_{s_k}(\omega) \to 0$. Because $S_{s_k}(\omega) \geq Y_{s_k}(\omega)$ and also $S_{s_k}(\omega) < 0$, it follows that $S_{s_k}(\omega) \to 0$. Since $S$ is nondecreasing, this implies $S_t(\omega) \uparrow 0$ as $t \uparrow \tau(\omega)$. Therefore, for every sufficiently large $n$, the level $-1/n$ is reached at some time strictly before $\tau(\omega)$, so $\sigma_n(\omega) < \tau(\omega)$. For such $n$, and in particular once $n > \tau(\omega)$, we have $\rho_n(\omega) = \sigma_n(\omega) < \tau(\omega) = \tau_L(\omega)$. As the levels $-1/n$ increase to 0, the corresponding hitting times increase to $\tau(\omega)$. Hence, $\rho_n(\omega) \uparrow \tau(\omega) = \tau_L(\omega)$.

If $\omega \notin L$ and $\tau(\omega) < \infty$, then $\tau_L(\omega) = \infty$. Since $\omega \notin L$ but $\tau(\omega) < \infty$, one has $Y_{\tau-}(\omega) < 0$. Choose $\varepsilon > 0$ such that $Y_{\tau-}(\omega) \leq -2\varepsilon$. By the definition of the left limit, there exists a left-neighbourhood of $\tau(\omega)$ on which $Y_t(\omega) \leq -\varepsilon$. In particular, near $\tau(\omega)$ the path remains uniformly below 0. Combined with Assumption 2.5, this shows that for all sufficiently large $n$, the level $-1/n$ cannot be reached before $\tau(\omega)$. Thus, $\sigma_n(\omega) \geq \tau(\omega)$ for all sufficiently large $n$, and therefore $\rho_n(\omega) \to \infty = \tau_L(\omega)$.

If $\tau = \infty$, then by Assumption 2.5, $S_t(\omega) < 0$ for every finite $t$. Fix $T < \infty$. Since $S_T(\omega) < 0$, one can choose $n$ so large that $-\frac{1}{n} > S_T(\omega)$. Then, by definition of $\sigma_n$, $\sigma_n(\omega) > T$. For all sufficiently large $n$, also $n > T$, so $\rho_n(\omega) = \sigma_n(\omega) \wedge n > T$. Since $T$ was arbitrary, it follows that $\rho_n(\omega) \to \infty = \tau_L(\omega)$.

Combining the three cases, we conclude that $\rho_n < \tau_L$ on $\{\tau_L < \infty\}$, $\rho_n \uparrow \tau_L$ pointwise. Therefore $\tau_L$ is predictable. □

*Corollary 2.11* (Sharpness of Assumption 2.5 for the canonical announcing sequence). *The sequence $(\rho_n)$ announces $\tau_L$ if and only if Assumption 2.5 holds.*

*Proof.* Sufficiency is exactly Theorem 2.10. Conversely, suppose that Assumption 2.5 fails. Then, there exist a sample path $w$ and a finite time $t$ with $t < \tau(w)$ such that $\sup_{0 \leq s \leq t} Y_s(w) = 0$. For every $n \geq 1$ we therefore have $S_t(w) \geq 0 \geq -1/n$, so by definition $\sigma_n(w) \leq t$. Hence $\sup_n \sigma_n(w) \leq t < \tau(w)$. If $w$ belongs to $L$, then $\tau_L(w) = \tau(w) > t$; if $w$ does not belong to $L$, then $\tau_L(w) = \infty$. In either case the sequence $(\sigma_n(w))$ cannot converge to $\tau_L(w)$. Therefore the canonical running-supremum sequence announces $\tau_L$ only if Assumption 2.5 holds. □

*Remark 2.12* Assumption 2.5 is exactly the condition needed for the running-supremum approximation to isolate the left-contact component. Without it, earlier left contacts followed by retreat may cause the times $\sigma_n$ to miss the true first-passage time.

## 2.4 The totally inaccessible side

**Proposition 2.13** (Gap from the left forces a jump). *On the event G, one has $\Delta X_\tau > 0$.*

*Proof.* On $G$, $Y_{\tau-} < 0 \leq Y_\tau$, hence $\Delta Y_\tau = Y_\tau - Y_{\tau-} > 0$. Since $b$ is continuous, $\Delta Y_\tau = \Delta X_\tau$. □

**Proposition 2.14** (What a predictable time would have to look like). *Let $\sigma$ be a predictable stopping time. Then,*

$$G \cap \{\tau = \sigma < \infty\} \subseteq \{\Delta X_\sigma > 0\}, \qquad L \cap \{\tau = \sigma < \infty\} \subseteq \{X_{\sigma-} = b(\sigma)\}. \tag{18}$$

*Proof.* On $G \cap \{\tau = \sigma < \infty\}$, Proposition 2.13 gives $\Delta X_\sigma > 0$. On $L \cap \{\tau = \sigma < \infty\}$, $Y_{\sigma-} = 0$, hence $X_{\sigma-} = b(\sigma)$. □

**Proposition 2.15** (Structural exclusion of predictable gap-crossings). *Assume that for every predictable stopping time $\sigma$,*

$$\mathbb{P}(\{\sigma < \infty\} \cap \{X_{\sigma-} < b(\sigma) \leq X_\sigma\}) = 0. \tag{19}$$

Then, for every predictable stopping time $\sigma$,

$$\mathbb{P}(\{\tau = \sigma < \infty\} \cap G) = 0. \tag{20}$$

*Proof.* On $\{\tau = \sigma < \infty\} \cap G$, one has $X_{\sigma-} < b(\sigma) \leq X_\sigma$ by definition of $G$. The hypothesis gives the desired zero-probability conclusion. □

**Corollary 2.16** (No jumps at predictable times). *If $\mathbb{P}(\Delta X_\sigma \neq 0, \sigma < \infty) = 0$ for every predictable stopping time $\sigma$, then*

$$\mathbb{P}(\{\tau = \sigma < \infty\} \cap G) = 0 \quad \text{for every predictable stopping time } \sigma. \tag{21}$$

*In particular, the gap-from-the-left mode is totally inaccessible in the natural eventwise sense.*

*Proof.* On $G \cap \{\tau = \sigma < \infty\}$, Proposition 2.13 yields $\Delta X_\sigma > 0$, contradicting the hypothesis. □

Now, define the restricted gap time by:

$$\tau_G(\omega) := \begin{cases} \tau(\omega), & \omega \in G, \\ \infty, & \omega \notin G. \end{cases}$$

**Proposition 2.17** (Stopping-time property of $\tau_G$). *The restricted random time $\tau_G$ is a stopping time.*

*Proof.* Since $L \in \mathcal{F}_\tau$ and

$$G = \{\tau < \infty\} \setminus L,$$

we have $G \in \mathcal{F}_\tau$. Hence, for every $t \geq 0$,

$$\{\tau_G \leq t\} = G \cap \{\tau \leq t\} \in \mathcal{F}_t.$$

Therefore, $\tau_G$ is a stopping time. □

**Corollary 2.18** (Classical total inaccessibility of the gap component). *Under the hypothesis of Proposition 2.15, the restricted random time $\tau_G$ is a totally inaccessible stopping time. In particular, this holds if $X$ has no jumps at predictable times.*

*Proof.* Let $\sigma$ be a predictable stopping time. Since $\tau_G$ is a stopping time,

$$\{\tau_G = \sigma < \infty\} = G \cap \{\tau = \sigma < \infty\}.$$

By Proposition 2.15, and in particular by Corollary 2.16 under the no-predictable-jumps hypothesis,

$$\mathbb{P}(G \cap \{\tau = \sigma < \infty\}) = 0.$$

Hence

$$\mathbb{P}(\tau_G = \sigma < \infty) = 0$$

for every predictable stopping time $\sigma$. Therefore $\tau_G$ is totally inaccessible. □

## 2.5 Compensator decomposition of the default indicator

**Theorem 2.19** (Compensator decomposition of the default indicator). *Define*

$$A_t := \mathbf{1}_{\{\tau \leq t\}}, \qquad A_t^L := \mathbf{1}_{\{\tau_L \leq t\}}, \qquad A_t^G := \mathbf{1}_{\{\tau_G \leq t\}}.$$

*Then $A = A^L + A^G$. Let $\Lambda, \Lambda^L, \Lambda^G$ denote the dual predictable projections of $A, A^L, A^G$, respectively. Then*

$$\Lambda = \Lambda^L + \Lambda^G.$$

*Moreover, $\Lambda^L = A^L$, so $\Lambda^L$ is purely discontinuous and has a single unit jump at $\tau_L$ on $L$, whereas $\Lambda^G$ is continuous.*

*Proof.* Since exactly one of $\tau_L$ and $\tau_G$ is finite on each path and, when finite, equals $\tau$, one has

$$A = A^L + A^G.$$

By linearity of dual predictable projection,

$$\Lambda = \Lambda^L + \Lambda^G.$$

Now $\tau_L$ is predictable, so the process $A^L$ is itself predictable. Therefore its dual predictable projection is just itself:

$$\Lambda^L = A^L.$$

In particular, $\Lambda^L$ is purely discontinuous, with a single jump of size 1 at $\tau_L$ on $L$.

Next consider $A^G$. By construction, $\tau_G$ is a totally inaccessible stopping time. Let

$$M^G := A^G - \Lambda^G.$$

Then $M^G$ is a local martingale. We claim that $\Lambda^G$ is continuous.

Let $\sigma$ be any predictable stopping time. Since

$$\Delta A^G_\sigma = \mathbf{1}_{\{\tau_G = \sigma < \infty\}},$$

and $\tau_G$ is totally inaccessible, one has

$$\mathbb{P}(\tau_G = \sigma < \infty) = 0,$$

hence

$$\Delta A^G_\sigma = 0 \quad \text{a.s.}$$

Because $\Lambda^G$ is predictable and $M^G = A^G - \Lambda^G$ is a local martingale, the predictable jump of $\Lambda^G$ at $\sigma$ is given by

$$\Delta \Lambda^G_\sigma = \mathbb{E}[\Delta A^G_\sigma \mid \mathcal{F}_{\sigma-}].$$

Therefore

$$\Delta \Lambda^G_\sigma = 0 \quad \text{a.s.}$$

for every predictable stopping time $\sigma$.

Now, $\Lambda^G$ is a predictable càdlàg increasing process. Any jump of such a process occurs at a predictable stopping time. Since we have shown that all predictable jumps vanish, $\Lambda^G$ has no jumps at all. Hence $\Lambda^G$ is continuous. □

# 3 From the fourfold decomposition to the contact-versus-overshoot dichotomy

## 3.1 Recovering trichotomy and dichotomy

**Assumption 3.1** (No jump from the graph at first passage). On $L$, $Y_\tau = 0$. Equivalently, $C_+ = \varnothing$.

**Proposition 3.2** (Trichotomy). *Under Assumption 3.1, $\{\tau < \infty\} = C_0 \dot\cup J_0 \dot\cup J_+$.*

*Proof.* By Proposition 2.2, $\{\tau < \infty\} = C_0 \dot\cup C_+ \dot\cup J_0 \dot\cup J_+$. Under Assumption 3.1, $C_+ = \varnothing$. □

**Assumption 3.3** (No exact hit by jump). $\mathbb{P}(J_0) = 0$.

**Corollary 3.4** (Contact versus overshoot dichotomy). *Under Assumptions 3.1 and 3.3, $\{\tau < \infty\} = C_0 \dot\cup J_+$ a.s.*

*Proof.* The trichotomy gives $\{\tau < \infty\} = C_0 \dot\cup J_0 \dot\cup J_+$. Assumption 3.3 removes $J_0$. □

**Proposition 3.5** (Compensator criterion for no exact hit). *Assume that $X$ is a semimartingale with jump measure $\mu^X$ and predictable compensator $\nu^X$. Define the predictable set*

$$\Gamma_0 := \{(\omega, t, x) : X_{t-}(\omega) < b(t),\ x = b(t) - X_{t-}(\omega)\}. \tag{22}$$

*If $\mathbf{1}_{\Gamma_0} * \nu^X \equiv 0$, then $\mathbb{P}(J_0) = 0$. In particular, this holds whenever $\nu^X$ admits a predictable disintegration $\nu^X(\omega, dt, dx) = K(\omega, t, dx)\, dA_t(\omega)$ with $K(\omega, t, \cdot)$ diffuse for $\mathbb{P} \otimes dA$-almost every $(\omega, t)$.*

*Proof.* Since $X_-$ is predictable and $b$ is continuous deterministic, $\Gamma_0$ is predictable. The process $N_t := \left(\mathbf{1}_{\Gamma_0} * \mu^X\right)_t$ counts jumps of $X$ landing exactly on the barrier from below. Under $\mathbf{1}_{\Gamma_0} * \nu^X \equiv 0$, $N$ is a nonnegative local martingale of finite variation, hence identically zero. On $J_0$, one has $\Delta X_\tau = b(\tau) - X_{\tau-}$, so $(\tau, \Delta X_\tau) \in \Gamma_0$. In particular, the jump at time $\tau$ contributes one point to the counting process $N$, so $N_\tau \geq 1$. But this contradicts the fact that $N \equiv 0$. Hence the event $J_0$ is impossible, and therefore $\mathbb{P}(J_0) = 0$. □

## 3.2 Characteristic criteria for Assumption 2.5 and for the jump modes

We now strengthen the preceding discussion in the semimartingale setting.

**Proposition 3.6** (Sharp pathwise reformulation of Assumption 2.5). *The following are equivalent:*

1. *Assumption 2.5 holds;*
2. $Y_{u-} < 0$ *for every finite* $u < \tau$;
3. *there is no finite* $u < \tau$ *such that* $Y_{u-} = 0 > Y_u$.

*Equivalently, Assumption 2.5 holds if and only if, before the first-passage time, the path never makes a downward jump away from the graph of the barrier.*

*Proof.* The implication (2) ⇒ (1) is Proposition 2.6. For (1) ⇒ (2), assume 2.5 and suppose $Y_{u-} = 0$ for some finite $u < \tau$. Since $u < \tau$, one has $Y_s < 0$ for all $s \leq u$, but $Y_{u-} = 0$ yields a sequence $s_n \uparrow u$ with $Y_{s_n} \to 0$. Hence $\sup_{0 \leq s \leq u} Y_s = 0$, contradicting 2.5 on $\{\tau > u\}$. The equivalence (2) ⇔ (3) is immediate because $u < \tau$ implies $Y_u < 0$. □

**Proposition 3.7** (Compensator criterion for Assumption 2.5). *Assume that X is a semimartingale. Define*

$$\Gamma^- := \{(\omega, t, x): t < \tau(\omega),\ X_{t-}(\omega) = b(t),\ x < 0\}. \tag{23}$$

*If* $\mathbf{1}_{\Gamma^-} * \nu^X \equiv 0$*, then Assumption 2.5 holds.*

*Proof.* Since $X_-$ is predictable, $b$ is continuous deterministic, and $\tau$ is a stopping time, the set $\Gamma^-$ is predictable. The counting process $N_t^- := (\mathbf{1}_{\Gamma^-} * \mu^X)_t$ counts negative jumps from the barrier graph before $\tau$. By assumption, is predictable compensator $\mathbf{1}_{\Gamma^-} * \nu^X \equiv 0$. Therefore $N^-$ is a nonnegative local martingale of finite variation. Hence it must be identically zero: $N^- \equiv 0$. We now prove Assumption 2.5 by contradiction. Suppose Assumption 2.5 fails. Then, by Proposition 3.6, there exists a finite time $u < \tau$ such that $Y_{u-} = 0 > Y_u$. Since $Y_t = X_t - b(t)$, this is equivalent to $X_{u-} = b(u)$, $\Delta X_u < 0$. Because $u < \tau$, the jump occurs before first passage. Therefore, $(u, \Delta X_u) \in \Gamma^-$. But then the counting process $N^-$ must record this jump, so in particular $N_u^- \geq 1$, contradicting the fact that $N^- \equiv 0$. □

**Corollary 3.8.** *If X has no negative jumps, then Assumption 2.5 holds.*

**Theorem 3.9** (Compensator formulas for the jump crossing modes). *Assume that X is a semimartingale. Define the predictable sets*

$$\begin{aligned}
\Gamma_+ &:= \{(\omega, t, x): X_{t-} = b(t),\ x > 0\}, \\
\Gamma_0 &:= \{(\omega, t, x): X_{t-} < b(t),\ x = b(t) - X_{t-}\}, \\
\Gamma_> &:= \{(\omega, t, x): X_{t-} < b(t),\ x > b(t) - X_{t-}\}.
\end{aligned} \tag{24}$$

*Then*

$$\mathbf{1}_{C_+} = \left(\mathbf{1}_{\Gamma_+} * \mu^X\right)_\tau, \qquad \mathbf{1}_{J_0} = \left(\mathbf{1}_{\Gamma_0} * \mu^X\right)_\tau, \qquad \mathbf{1}_{J_+} = \left(\mathbf{1}_{\Gamma_>} * \mu^X\right)_\tau. \tag{25}$$

*Consequently,*

$$\mathbb{P}(C_+) = \mathbb{E}\left[\left(\mathbf{1}_{\Gamma_+} * \nu^X\right)_\tau\right], \quad \mathbb{P}(J_0) = \mathbb{E}\left[\left(\mathbf{1}_{\Gamma_0} * \nu^X\right)_\tau\right], \quad \mathbb{P}(J_+) = \mathbb{E}\left[\left(\mathbf{1}_{\Gamma_>} * \nu^X\right)_\tau\right]. \tag{26}$$

*Equivalently,*

$$\mathbb{P}(C_+) = \mathbb{E}\left[\int_0^\tau \mathbf{1}_{\{X_{t-}=b(t)\}}\, \nu^X(dt,(0,\infty))\right], \tag{27}$$

$$\mathbb{P}(J_0) = \mathbb{E}\left[\int_0^\tau \mathbf{1}_{\{X_{t-}<b(t)\}}\, \nu^X(dt,\{b(t)-X_{t-}\})\right], \tag{28}$$

$$\mathbb{P}(J_+) = \mathbb{E}\left[\int_0^\tau \mathbf{1}_{\{X_{t-}<b(t)\}}\, \nu^X(dt,(b(t)-X_{t-},\infty))\right]. \tag{29}$$

*Proof.* We prove the identity for $C_+$; the arguments for $J_0$ and $J_+$ are identical. The process $\left(\mathbf{1}_{\Gamma_+} * \mu^X\right)_t$ counts the number of jumps of $X$ up to time $t$ such that just before the jump the process lies on the barrier graph and the jump is strictly positive. Suppose first that the event $C_+$ occurs. Then, by definition of $C_+$, $X_{\tau-} = b(\tau)$, $\Delta X_\tau > 0$. Hence the jump of $X$ at time $\tau$ belongs to $\Gamma_+$, so $\left(\mathbf{1}_{\Gamma_+} * \mu^X\right)_\tau \geq 1$. On the other hand, there can be no earlier time $s < \tau$ such that $X_{s-} = b(s)$, $\Delta X_s > 0$, because then $X_s = X_{s-} + \Delta X_s > b(s)$, so first passage would already have occurred at time $s$, contradicting $s < \tau$. Therefore the only possible jump counted by $\left(\mathbf{1}_{\Gamma_+} * \mu^X\right)_\tau$ is the jump at time $\tau$, and thus $\left(\mathbf{1}_{\Gamma_+} * \mu^X\right)_\tau = 1$. This proves that

$$\mathbf{1}_{C_+} \leq \left(\mathbf{1}_{\Gamma_+} * \mu^X\right)_\tau.$$

Conversely, suppose that $\left(\mathbf{1}_{\Gamma_+} * \mu^X\right)_\tau = 1$. Then there exists a jump time $s \leq \tau$ such that $X_{s-} = b(s)$, $\Delta X_s > 0$. If $s < \tau$, then $X_s > b(s)$, which again contradicts the definition of $\tau$. Hence necessarily $s = \tau$, and therefore $X_{\tau-} = b(\tau)$, $\Delta X_\tau > 0$. That is exactly the event $C_+$. So,

$$\left(\mathbf{1}_{\Gamma_+} * \mu^X\right)_\tau \leq \mathbf{1}_{C_+}.$$

The arguments for $J_0$ and $J_+$ are identical. Taking expectations and applying the compensation formula yields (26)–(29). $\square$

**Corollary 3.10** (Characteristic conditions for the classical dichotomy). *Assume that $X$ is a semimartingale and that:*

1. *on $\{t < \tau,\, X_{t-} = b(t)\}$, the compensator assigns no mass to $(-\infty,0)$;*

2. on $\{X_{t-} = b(t)\}$, the compensator assigns no mass to $(0, \infty)$;

3. on $\{X_{t-} < b(t)\}$, the compensator assigns no mass to $\{b(t) - X_{t-}\}$.

Then Assumption 2.5 holds, $\mathbb{P}(C_+) = 0$, $\mathbb{P}(J_0) = 0$, and $\{\tau < \infty\} = C_0 \dot\cup J_+$ a.s. If additionally X has no jumps at predictable times, then $\tau_L$ is predictable and the $J_+$-component is totally inaccessible in the eventwise sense.

## 3.3 Examples

**Example 3.11** (Continuous semimartingale). Let X be continuous. Then $Y_{u-} = Y_u < 0$ for every $u < \tau$, so Assumption 2.5 holds. Continuity excludes jumps at $\tau$, hence $C_+ = J_0 = J_+ = \emptyset$ and $\{\tau < \infty\} = C_0$: first passage is entirely by continuous contact. This is exactly the classical notion of creeping.

**Example 3.12** (Compound Poisson process). Let $X_t = \sum_{k=1}^{N_t} U_k$, where N has intensity $\lambda > 0$ and the $U_k$ are i.i.d., independent of N. Let $b(t) := a > 0$. There is no continuous motion at all: the process just sits still between jump times. This is the pure-jump analogue of the continuous-semimartingale example 3.11, where instead all the jump-driven modes vanish. Sample paths look like upward step functions: flat, then jump, then flat, then jump, and so on. The event $C_+$ cannot happen because for $Y_{\tau-} = 0$, we would need $X_{\tau-} = a$; but if the process is already exactly at the barrier just before $\tau$, then since the path is constant between jumps, it must already have hit the barrier earlier. Thus, $C_0 = C_+ = \emptyset$ and $\{\tau < \infty\} = J_0 \dot\cup J_+$. The compensator is $\nu^X(dt, dx) = \lambda F(dx)\, dt$, where F is the law of $U_1$. Then

$$\mathbb{P}(J_0) = \mathbb{E}\left[\int_0^\tau \lambda\, F(\{a - X_{t-}\})\, dt\right], \qquad \mathbb{P}(J_+) = \mathbb{E}\left[\int_0^\tau \lambda\, F((a - X_{t-}, \infty))\, dt\right]. \tag{30}$$

Suppose the jump-size distribution is lattice-valued, for example $U_k \in \{1, 2, 3, \ldots\}$ or more generally takes values in some arithmetic progression. Then the singleton probability $F(\{a - X_{t-}\})$ need not be zero. If the current gap $a - X_{t-}$ happens to be one of the lattice values, there may be a strictly positive chance that the next jump fills the gap exactly. So exact hit by jump can genuinely occur with positive probability: $\mathbb{P}(J_0) > 0$ is possible. If F is diffuse, then every singleton has zero mass. In particular, $F(\{a - X_{t-}\}) = 0$ for every t. Hence the only remaining possibility is strict overshoot: $\{\tau < \infty\} = J_+$ a.s. The conclusion is clean and intuitive: for a pure-jump process with diffuse jump sizes, first passage cannot occur by landing exactly on the barrier; it must occur by overshoot.

**Example 3.13** (A path in the $C_+$ mode). Let $b(t) \equiv 0$ and $X_t = -1$ for $t \in [0, 1)$, $X_t = t - 2$ for $t \in [1, 2)$, and $X_t = 1$ for $t \geq 2$. Then $\tau = 2$, $Y_{\tau-} = 0$, $Y_\tau = 1 > 0$. The path arrives at the barrier continuously, but crosses by jumping. So the path belongs to $C_+$. This shows that the three-way split into $C_0, J_0, J_+$ is not exhaustive and that $C_+$ is a genuine fourth mode that cannot be ignored without the extra assumption 3.1.

**Example 3.14** (Failure of Assumption 2.5). Let $b(t) \equiv 0$ and define $X_t = t - 1$ for $t \in [0,1)$, $X_1 = -\frac{1}{2}$, $X_t = t - 2$ for $t \in (1,2)$, and $X_t = 0$ for $t \geq 2$. Then $X_{1-} = 0 > X_1 = -\frac{1}{2}$: the path makes a premature left contact at $t = 1$ and then retreats. The actual first passage is at $\tau = 2$. The running supremum satisfies $S_t = 0$ for all $t \geq 1$, so $\sigma_n \leq 1 < 2 = \tau$ for all $n$: the announcing sequence fails to converge to $\tau$. This example shows that Assumption 2.5 is needed because the running supremum remembers the earlier contact and the times $\sigma_n$ "lock onto" the wrong moment.

**Example 3.15** (Spectrally negative tempered stable process). Let $X$ be a spectrally negative tempered stable Lévy process with constant barrier $b(t) \equiv a > 0$. Since the process has no positive jumps, none of the jump-driven upward crossing modes can occur. Indeed: the event $C_+$ would require a positive jump from the barrier graph; the event $J_0$ would require a positive jump from below landing exactly on the barrier; the event $J_+$ would require a positive jump from below overshooting the barrier. The compensator formulas give $\mathbb{P}(C_+) = \mathbb{P}(J_0) = \mathbb{P}(J_+) = 0$. Hence, the fourfold decomposition reduces to $\{\tau < \infty\} = C_0$ a.s.: upward first passage, when it occurs, is entirely by creeping. The compensator formulas recover this classical fact, but in the stronger form that the three jump-driven modes vanish separately.

## 4 Affine jump-diffusion model and general overshoot equation

This section develops an explicit model in which the compensator formulas of Section 3 lead to a genuinely nontrivial analytic problem. We consider an affine jump-diffusion with upward exponential jumps. The memoryless jump law still yields an explicit overshoot law, while the state-dependent drift produces a genuinely non-homogeneous boundary-value problem.

### 4.1 Model specification

Let

$$X_t = x + \int_0^t (\alpha + \beta X_s) \, ds + \sigma W_t + \sum_{k=1}^{N_t} U_k, \qquad t \geq 0, \tag{31}$$

where $x < a$, $\alpha, \beta \in \mathbb{R}$, $\sigma > 0$, $W$ is a standard Brownian motion, $N$ is a Poisson process of intensity $\lambda > 0$ independent of $W$, and $(U_k)_{k \geq 1}$ are i.i.d. with $U_1 \sim \text{Exp}(\eta)$, $\eta > 0$, independent of $(W, N)$. The barrier is constant, $b(t) \equiv a$, and we write

$$\tau_a := \inf\{t \geq 0 \colon X_t \geq a\}.$$

The jump compensator is

$$\nu^X(dt, dz) = \lambda \eta \, e^{-\eta z} \, \mathbf{1}_{\{z > 0\}} \, dz \, dt. \tag{32}$$

Since the jumps are positive and the kernel is diffuse, Corollary 3.8 gives Assumption 2.5, Proposition 3.5 gives $\mathbb{P}(J_0) = 0$, and (27) gives $\mathbb{P}(C_+) = 0$ because the occupation time at a fixed level vanishes for a semimartingale with a nondegenerate continuous martingale part. Hence

$$\{\tau_a < \infty\} = C_0 \,\dot\cup\, J_+ \quad \text{a.s.} \tag{33}$$

For $q \geq 0$, define the discounted mode-separated transforms

$$G_q(x) := \mathbb{E}_x\!\left[e^{-q\tau_a}\mathbf{1}_{J_+}\right], \qquad H_q(x) := \mathbb{E}_x\!\left[e^{-q\tau_a}\mathbf{1}_{\tau_a<\infty}\right], \tag{34}$$

$$F_q(x) := H_q(x) - G_q(x).$$

Thus $F_q(x) = \mathbb{E}_x\!\left[e^{-q\tau_a}\mathbf{1}_{C_0}\right]$ by (33).

## 4.2 The overshoot distribution

**Theorem 4.1** (Overshoot distribution on $J_+$). *Under the model (31), conditional on $J_+$, the overshoot $X_{\tau_a} - a$ is $\mathrm{Exp}(\eta)$-distributed and is independent of the pre-jump level $X_{\tau_a-}$. More precisely, for every measurable $g\colon[0,\infty) \to [0,\infty)$,*

$$\mathbb{E}_x\!\left[g(X_{\tau_a} - a)\,\mathbf{1}_{J_+}\right] = \mathbb{P}_x(J_+) \cdot \int_0^\infty g(u)\,\eta\,e^{-\eta u}\,du. \tag{35}$$

*Proof.* By Theorem 3.9, $\mathbf{1}_{J_+} = \left(\mathbf{1}_{\Gamma_>} * \mu^X\right)_{\tau_a}$. Hence, for measurable $g \geq 0$,

$$\mathbb{E}_x\!\left[g(X_{\tau_a} - a)\,\mathbf{1}_{J_+}\right] = \mathbb{E}_x\!\left[\int_0^{\tau_a}\!\int_0^\infty g(X_{t-} + z - a)\,\mathbf{1}_{\{X_{t-}<a,\,z>a-X_{t-}\}}\,\nu^X(dt, dz)\right]. \tag{36}$$

Substituting (32) and setting $u = z - (a - X_{t-})$, the inner integral becomes

$$\lambda\,e^{-\eta(a-X_{t-})}\int_0^\infty g(u)\,\eta\,e^{-\eta u}\,du. \tag{37}$$

The exponential kernel factorizes the $g$-dependence from the level-dependence. Substituting back,

$$\mathbb{E}_x\!\left[g(X_{\tau_a} - a)\,\mathbf{1}_{J_+}\right] = \left(\int_0^\infty g(u)\,\eta\,e^{-\eta u}\,du\right)\mathbb{E}_x\!\left[\int_0^{\tau_a}\lambda\,e^{-\eta(a-X_{t-})}\,\mathbf{1}_{\{X_{t-}<a\}}\,dt\right].$$

Setting $g \equiv 1$ identifies the second factor as $\mathbb{P}_x(J_+)$. □

*Remark 4.2.* The factorization (37) is specific to the exponential distribution. For a general positive jump-size density $f$, the same compensation argument gives

$$\mathbb{E}_x[g(X_{\tau_a} - a)\mathbf{1}_{J_+}] = \mathbb{E}_x\left[\int_0^{\tau_a} \lambda \mathbf{1}_{\{X_{t-}<a\}} \int_0^\infty g(u) f(u + a - X_{t-})\, du\, dt\right],$$

so the overshoot law on $J_+$ typically depends on $X_{\tau_a-}$.

## 4.3 Integro-differential equation and ODE reduction

**Proposition 4.3** (Integro-differential equation). *For every $q \geq 0$, the function $G_q$ satisfies, for $x < a$,*

$$\frac{1}{2}\sigma^2 G_q''(x) + (\alpha + \beta x) G_q'(x) + \lambda \int_0^{a-x} G_q(x+z)\, \eta\, e^{-\eta z}\, dz - (\lambda + q) G_q(x) \qquad (38)$$
$$+ \lambda e^{-\eta(a-x)} = 0.$$

*Equivalently, if*

$$I_q(x) := \int_x^a G_q(y)\, \eta\, e^{-\eta(y-x)}\, dy, \qquad (39)$$

*then*

$$\frac{1}{2}\sigma^2 G_q''(x) + (\alpha + \beta x) G_q'(x) + \lambda I_q(x) - (\lambda + q) G_q(x) + \lambda e^{-\eta(a-x)} = 0, \qquad (40)$$

*with*

$$I_q'(x) = \eta I_q(x) - \eta G_q(x). \qquad (41)$$

*Proof.*

Before the hitting time $\tau_a$, the process evolves as the affine jump-diffusion

$$dX_t = (\alpha + \beta X_t)\, dt + \sigma\, dW_t + dJ_t,$$

but is *killed* as soon as it reaches or exceeds the level $a$. Consequently, when the process is currently at a point $x < a$, jumps of size $z < a - x$ keep the process alive, whereas jumps of size $z > a - x$ terminate the evolution and realize the event $J_+$ immediately.

It is therefore convenient to introduce the killed generator

$$\mathcal{L}^a f(x) := (\alpha + \beta x) f'(x) + \frac{1}{2}\sigma^2 f''(x) + \lambda \int_0^{a-x} f(x+z)\, \eta e^{-\eta z}\, dz - \lambda f(x), \qquad x < a,$$

acting on sufficiently smooth test functions $f$. The truncation at $a - x$ reflects the fact that jumps above this threshold no longer contribute to the continuation of the process: they kill it.

Now, when the process is at $x < a$, an overshooting jump occurs at instantaneous rate
$$r(x) := \lambda\, \mathbb{P}(U > a - x) = \lambda e^{-\eta(a-x)},$$
because the jump sizes are exponentially distributed with parameter $\eta$. Such a jump leads to immediate first passage and contributes the payoff 1 in the definition of $G_q$. Thus $G_q$ should satisfy the resolvent equation for the killed process with killing rate $q$ and source term $r$, namely
$$(\mathcal{L}^a - q)G_q(x) + r(x) = 0, \qquad x < a.$$
Expanding this identity yields (38). The reformulation (40) follows from the change of variables
$$\int_0^{a-x} G_q(x+z)\, \eta e^{-\eta z}\, dz = \int_x^a G_q(y)\, \eta e^{-\eta(y-x)}\, dy = I_q(x).$$
Differentiating (39) gives (41). $\square$

**Lemma 4.4** (Correct ODE reduction and compatibility). *Define the residual*

$$R_q(x) := \frac{1}{2}\sigma^2 G_q''(x) + (\alpha + \beta x)G_q'(x) + \lambda I_q(x) - (\lambda + q)G_q(x) + \lambda e^{-\eta(a-x)}. \tag{42}$$

*Then*

$$R_q'^{(x)} - \eta R_q(x) = \frac{1}{2}\sigma^2 G_q'''(x) \tag{43}$$

$$+ \left(\alpha + \beta x - \frac{1}{2}\eta\sigma^2\right) G_q''(x) \quad + (\beta - \eta(\alpha + \beta x) - \lambda - q)G_q'(x) + \eta q\, G_q(x).$$

*Consequently, the integro-differential equation (38) is equivalent to the third-order ODE*

$$\frac{1}{2}\sigma^2 G_q'''(x) + \left(\alpha + \beta x - \frac{1}{2}\eta\sigma^2\right) G_q''(x) \tag{44}$$

$$+ (\beta - \eta(\alpha + \beta x) - \lambda - q)G_q'(x) + \eta q\, G_q(x) = 0, \qquad x < a,$$

*together with one compatibility condition. Since $I_q(a^-) = 0$ and $G_q(a^-) = 0$, this may be taken as*

$$\frac{1}{2}\sigma^2 G_q''(a^-) + (\alpha + \beta a)G_q'(a^-) + \lambda = 0. \tag{45}$$

*Proof.* By construction, the integro-differential equation is exactly the statement

$$R_q(x) = 0, \qquad x < a.$$

We first compute $R_q'(x) - \eta R_q(x)$. Differentiating $R_q$ term by term, and using

$$\frac{d}{dx}\left((\alpha + \beta x)G_q'(x)\right) = (\alpha + \beta x)G_q''(x) + \beta G_q'(x),$$

$$I_q'(x) = \eta I_q(x) - \eta G_q(x),$$

we obtain:

$$R_q'(x) = \frac{1}{2}\sigma^2 G_q'''(x) + (\alpha + \beta x)G_q''(x) + \beta G_q'(x) + \lambda \eta I_q(x) - \lambda \eta G_q(x) - (\lambda + q)G_q'(x)$$
$$+ \eta \lambda e^{-\eta(a-x)}.$$

Subtracting $\eta R_q(x)$, the terms involving $I_q$ and the exponential term cancel, and we obtain (43). Now suppose first that the integro-differential equation holds. Then, $R_q(x) = 0$ for all $x < a$, and therefore

$$R_q'(x) - \eta R_q(x) = 0.$$

Hence $G_q$ satisfies the third-order ODE in (44).

Conversely, suppose that $G_q$ satisfies the third-order ODE. Then the identity just proved shows that $R_q'(x) - \eta R_q(x) = 0$. This is a first-order linear differential equation for $R_q$, whose general solution is $R_q(x) = Ce^{\eta x}$ for some constant $C$. Thus the third-order ODE alone determines the residual only up to a multiplicative constant. To recover the original integro-differential equation, one must therefore impose one additional scalar condition.

Since

$$I_q(a^-) = \int_a^a G_q(y)\,\eta e^{-\eta(y-a)}\,dy = 0,$$

and

$$G_q(a^-) = 0,$$

we have

$$R_q(a^-) = \frac{1}{2}\sigma^2 G_q''(a^-) + (\alpha + \beta a)G_q'(a^-) + \lambda.$$

Hence the boundary condition

$$\frac{1}{2}\sigma^2 G_q''(a^-) + (\alpha + \beta a)G_q'(a^-) + \lambda = 0$$

is exactly the statement

$$R_q(a^-) = 0.$$

But if $R_q(x) = Ce^{\eta x}$, then $R_q(a^-) = Ce^{\eta a}$. Therefore $R_q(a^-) = 0$ forces $C = 0$, and so

$$R_q(x) \equiv 0.$$

This is precisely the original integro-differential equation.

We have therefore shown that the integro-differential equation is equivalent to the third-order ODE together with one compatibility condition, which may be taken to be (45). □

*Remark 4.5* The same calculation applies to the total first-passage transform

$$H_q(x) = \mathbb{E}_x\big[e^{-q\tau_a}\mathbf{1}_{\{\tau_a < \infty\}}\big].$$

The resulting third-order ODE is again (44), but with boundary data

$$H_q(a^-) = 1, \qquad \frac{1}{2}\sigma^2 H_q''(a^-) + (\alpha + \beta a)H_q'(a^-) - q = 0.$$

Hence $F_q = H_q - G_q$ gives the discounted left-contact component.

**Theorem 4.6** (Verification and uniqueness for the discounted affine OIDE). *Fix $q > 0$. Let $u \in C^2((-\infty, a)) \cap C((-\infty, a])$ be bounded, assume $u(a^-) = 0$, and suppose that for every $x < a$,*

$$\frac{1}{2}\sigma^2 u''(x) + (\alpha + \beta x)u'(x) + \lambda \int_0^{a-x} u(x+z)\eta e^{-\eta z}\,dz - (\lambda + q)u(x) + \lambda e^{-\eta(a-x)} \tag{46}$$
$$= 0.$$

*Then*

$$u(x) = \mathbb{E}_x\big[e^{-q\tau_a}\mathbf{1}_{J_+}\big] = G_q(x), \qquad x < a.$$

*In particular, for each $q > 0$, the discounted overshoot transform $G_q$ is the unique bounded classical solution of (46) with boundary condition $u(a^-) = 0$.*

*Proof.* Fix $q > 0$, and let $u$ satisfy the hypotheses stated above.

We first rewrite this equation in generator form. Define the killed generator

$$\mathcal{L}^a f(x) := (\alpha + \beta x)f'(x) + \frac{1}{2}\sigma^2 f''(x) + \lambda \int_0^{a-x} f(x+z)\eta e^{-\eta z}\,dz - \lambda f(x), \qquad x < a.$$

Then the assumed equation for $u$ is exactly

$$(\mathcal{L}^a - q)u(x) + \lambda e^{-\eta(a-x)} = 0, \qquad x < a.$$

Thus $u$ solves the discounted resolvent equation for the process killed at $\tau_a$, with running reward

$$r(x) := \lambda e^{-\eta(a-x)}.$$

Apply the discounted Dynkin formula to the killed process up to time $t \wedge \tau_a$. Since $u$ is a bounded classical solution, we obtain

$$u(x) = \mathbb{E}_x\left[e^{-q(t\wedge\tau_a)}u(X^a_{t\wedge\tau_a})\right] + \mathbb{E}_x\left[\int_0^{t\wedge\tau_a} e^{-q}\,r(X_{s-})\,ds\right].$$

Because $u$ is bounded and $q > 0$,

$$\left|e^{-q(t\wedge\tau_a)}u(X^a_{t\wedge\tau_a})\right| \leq \|u\|_\infty\, e^{-q(t\wedge\tau_a)},$$

and the right-hand side converges to 0 as $t \to \infty$. Hence, by dominated convergence,

$$u(x) = \mathbb{E}_x\left[\int_0^{\tau_a} e^{-qs}\,r(X_{s-})\,ds\right]$$

So every bounded classical solution of the OIDE admits the representation

$$u(x) = \mathbb{E}_x\left[\int_0^{\tau_a} e^{-qs}\,\lambda e^{-\eta(a-X_{s-})}\,ds\right]. \tag{47}$$

We now show that the same representation holds for $G_q(x) := \mathbb{E}_x\left[e^{-q\tau_a}\mathbf{1}_{J_+}\right]$.

On the event $J_+$, first passage is caused by a jump from a point strictly below the barrier, so

$$e^{-q\tau_a}\mathbf{1}_{J_+} = \sum_{t\leq\tau_a} e^{-qt}\,\mathbf{1}_{\{X_{t-}<a<X_{t-}+\Delta J_t\}}.$$

Since all jumps are upward and $\tau_a$ is the first crossing time, at most one jump can contribute to this sum, namely the jump occurring at $\tau_a$ on the overshoot event. Taking expectations and applying the compensation formula for the Poisson jump measure yields

$$\mathbb{E}_x\left[e^{-q\tau_a}\mathbf{1}_{J_+}\right] = \mathbb{E}_x\left[\int_0^{\tau_a} e^{-q}\,\lambda\,\mathbb{P}(U > a - X_{t-})\,dt\right].$$

Since $U \sim \text{Exp}(\eta)$, we have

$$G_q(x) = \mathbb{E}_x\left[\int_0^{\tau_a} e^{-qt}\,\lambda e^{-\eta(a-X_{t-})}\,dt\right].$$

Comparing this with the representation obtained for $u$, we conclude that

$$u(x) = G_q(x), \qquad x < a.$$

This proves the verification statement. The uniqueness claim is now immediate: if $u_1$ and $u_2$ are two bounded classical solutions of the OIDE with boundary condition $u(a^-) = 0$, then both equal $G_q$, hence $u_1 \equiv u_2$. □

## Section 5 — Discounted analysis in the mean-reverting case

Assume throughout this subsection that $\beta < 0$ and write $b := -\beta > 0$. For $q \geq 0$ let

$$w_q(x) := G_q'(x), \qquad v_q := -1 - \frac{\lambda + q}{b}, \tag{48}$$

and define the second-order differential operator

$$L_q := 1/2\,\sigma^2 \frac{d^2}{dx^2} + \left(\alpha - bx - \frac{1}{2}\eta\sigma^2\right)\frac{d}{dx} + (-b - \eta(\alpha - bx) - \lambda - q). \tag{49}$$

Since $G_q(a^-) = 0$, differentiating (44) and writing

$$G_q(x) = -\int_x^a w_q(y)\,dy$$

shows that $w_q$ satisfies

$$L_q w_q(x) = \eta q \int_x^a w_q(y)\,dy, \qquad x < a, \tag{50}$$

together with the boundary condition

$$\mathcal{B}_a w_q := \frac{1}{2}\sigma^2 w_q'(a^-) + (\alpha - ba) w_q(a^-) = -\lambda. \tag{51}$$

Thus the discounted problem is a Volterra perturbation of a local Weber-type equation.

Introduce the gauge transform and scaled variable

$$p(x) := -\frac{\alpha x}{\sigma^2} + \frac{bx^2}{2\sigma^2} + \frac{\eta x}{2}, \qquad z(x) := \frac{\sqrt{2}}{\sigma\sqrt{b}}\left(-bx + \alpha + \frac{1}{2}\eta\sigma^2\right). \tag{52}$$

Then the homogeneous equation $L_q u = 0$ is transformed by $u(x) = e^{p(x)}\phi(z(x))$ into

$$\phi_{zz}(z) + \left(v_q + \frac{1}{2} - \frac{z^2}{4}\right)\phi(z) = 0. \tag{53}$$

Accordingly, define

$$\psi_q(x) := e^{p(x)} D_{v_q}(z(x)), \qquad \tilde{\psi}_q(x) := e^{p(x)} D_{v_q}(-z(x)). \tag{54}$$

By the standard asymptotics of parabolic-cylinder functions, $\psi_q$ is the unique (up to scalar multiples) solution of $L_q u = 0$ that remains bounded as $x \to -\infty$.

**Lemma 5.1** (Boundary operator on the homogeneous basis). *For every $q \geq 0$,*

$$\mathcal{B}_a \psi_q = \frac{\sigma \sqrt{b}}{\sqrt{2}} e^{p(a)} D_{v_q+1}(z(a)), \tag{55}$$

*and*

$$\mathcal{B}_a \tilde{\psi}_q = -\frac{\sigma \sqrt{b}}{\sqrt{2}} e^{p(a)} D_{v_q+1}(-z(a)). \tag{56}$$

*Consequently, whenever $D_{v_q+1}(z(a)) \neq 0$, the function*

$$\chi_q(x) := \tilde{\psi}_q(x) + \frac{D_{v_q+1}(-z(a))}{D_{v_q+1}(z(a))} \psi_q(x) \tag{57}$$

*satisfies $L_q \chi_q = 0$ and $\mathcal{B}_a \chi_q = 0$.*

*Proof.* By construction, both $\psi_q$ and $\tilde{\psi}_q$ arise from the transformation of the homogeneous equation $L_q u = 0$ into the parabolic-cylinder equation. Hence

$$L_q \psi_q = 0, \qquad L_q \tilde{\psi}_q = 0.$$

It remains to compute the action of the boundary operator

$$\mathcal{B}_a w := \frac{1}{2} \sigma^2 w'(a^-) + (\alpha - ba) w(a^-)$$

on these two functions.

We begin with $\psi_q(x) = e^{p(x)} D_{v_q}(z(x))$. Differentiating and using the chain rule gives

$$\psi_q'(x) = e^{p(x)} \left( p'(x) D_{v_q}(z(x)) + D_{v_q}'(z(x)) z'(x) \right).$$

Therefore

$$\mathcal{B}_a \psi_q = e^{p(a)} \left[ \left( \frac{1}{2} \sigma^2 p'(a) + \alpha - ba \right) D_{v_q}(z(a)) + \frac{1}{2} \sigma^2 z'(a) D_{v_q}'(z(a)) \right].$$

Now, $z'(x) = -\frac{\sqrt{2b}}{\sigma}$, and from the definition of $p$ and $z$ one checks that

$$\frac{1}{2}\sigma^2 p'(a) + (\alpha - ba) = \frac{\sigma\sqrt{b}}{2\sqrt{2}} z(a).$$

Moreover, the standard derivative identity for parabolic-cylinder functions reads

$$D_\nu'(z) = \frac{z}{2} D_\nu(z) - D_{\nu+1}(z).$$

Substituting these identities into the preceding expression yields

$$\mathcal{B}_a \psi_q = e^{p(a)} \left[ \frac{\sigma\sqrt{b}}{2\sqrt{2}} z(a) D_{\nu_q}(z(a)) - \frac{\sigma\sqrt{2b}}{2} \left( \frac{z(a)}{2} D_{\nu_q}(z(a)) - D_{\nu_q+1}(z(a)) \right) \right].$$

The two terms involving $z(a) D_{\nu_q}(z(a))$ cancel, and we obtain the formula claimed in (55).

Next, consider $\tilde{\psi}_q(x) = e^{p(x)} D_{\nu_q}(-z(x))$. Differentiating gives

$$\tilde{\psi}_q'(x) = e^{p(x)} \left( p'(x) D_{\nu_q}(-z(x)) - D_{\nu_q}'(-z(x)) z'(x) \right),$$

because $\frac{d}{dx}(-z(x)) = -z'(x)$. Hence

$$\mathcal{B}_a \tilde{\psi}_q = e^{p(a)} \left[ \left( \frac{1}{2}\sigma^2 p'(a) + \alpha - ba \right) D_{\nu_q}(-z(a)) - \frac{1}{2}\sigma^2 z'(a) D_{\nu_q}'(-z(a)) \right].$$

Using again $\frac{1}{2}\sigma^2 p'(a) + (\alpha - ba) = \frac{\sigma\sqrt{b}}{2\sqrt{2}} z(a)$ and $z'(a) = -\frac{\sqrt{2b}}{\sigma}$, together with

$$D_{\nu_q}'(-z(a)) = -\frac{z(a)}{2} D_{\nu_q}(-z(a)) - D_{\nu_q+1}(-z(a)),$$

we obtain

$$\mathcal{B}_a \tilde{\psi}_q = e^{p(a)} \left[ \frac{\sigma\sqrt{b}}{2\sqrt{2}} z(a) D_{\nu_q}(-z(a)) + \frac{\sigma\sqrt{2b}}{2} \left( -\frac{z(a)}{2} D_{\nu_q}(-z(a)) - D_{\nu_q+1}(-z(a)) \right) \right].$$

Again the terms involving $z(a) D_{\nu_q}(-z(a))$ cancel, and we obtain the formula claimed in (56).

Finally, assume that $D_{\nu_q+1}(z(a)) \neq 0$, and define

$$\chi_q(x) := \tilde{\psi}_q(x) + \frac{D_{\nu_q+1}(-z(a))}{D_{\nu_q+1}(z(a))} \psi_q(x).$$

Since both $\psi_q$ and $\tilde{\psi}_q$ solve $L_q u = 0$, the same is true of $\chi_q$. Moreover, by linearity of $\mathcal{B}_a$,

$$\mathcal{B}_a \chi_q = \mathcal{B}_a \tilde{\psi}_q + \frac{D_{\nu_q+1}(-z(a))}{D_{\nu_q+1}(z(a))} \mathcal{B}_a \psi_q.$$

Substituting the two formulas just obtained, the two terms cancel exactly, so $\mathcal{B}_a \chi_q = 0$. □

**Theorem 5.2** (Green–Volterra representation for the discounted overshoot mode). *Assume $\beta < 0$, fix $q \geq 0$, and suppose that $D_{\nu_q+1}(z(a)) \neq 0$. Let*

$$W_q(x) := \psi_q(x)\chi_q'(x) - \psi_q'(x)\chi_q(x). \tag{58}$$

*Then the derivative $w_q = G_q'$ admits the representation*

$$w_q(x) = w_q^{(0)}(x) + \eta q \int_{-\infty}^{a} \Gamma_q(x, y) \left( \int_y^a w_q(u)\, du \right) dy, \qquad x < a, \tag{59}$$

*where*

$$w_q^{(0)}(x) = -\frac{\sqrt{2}\,\lambda}{\sigma\sqrt{b}} \exp(p(x) - p(a)) \frac{D_{\nu_q}(z(x))}{D_{\nu_q+1}(z(a))}, \tag{60}$$

*and the Green kernel is*

$$\Gamma_q(x, y) := \frac{2}{\sigma^2 W_q(y)} \begin{cases} \psi_q(x)\chi_q(y), & x \leq y < a, \\ \psi_q(y)\chi_q(x), & y \leq x < a. \end{cases} \tag{61}$$

*In particular, the discounted problem is reduced exactly to a second-order Weber equation together with the explicit Volterra perturbation appearing in (59).*

*Proof.* Assume $\beta < 0$, and write $b := -\beta > 0$. Recall from the preceding reduction that the derivative $w_q := G_q'$ satisfies

$$L_q w_q(x) = \eta q \int_x^a w_q(y)\, dy, \qquad x < a,$$

where

$$L_q = \frac{1}{2}\sigma^2 \frac{d^2}{dx^2} + \left( \alpha - bx - \frac{1}{2}\eta\sigma^2 \right) \frac{d}{dx} + (-b - \eta(\alpha - bx) - \lambda - q),$$

together with the Robin boundary condition

$$\mathcal{B}_a w_q = -\lambda, \qquad \mathcal{B}_a w := \frac{1}{2}\sigma^2 w'(a^-) + (\alpha - ba)w(a^-).$$

Set

$$f_q(x) := \eta q \int_x^a w_q(y)\, dy.$$

Then $w_q$ solves the linear inhomogeneous boundary-value problem

$$L_q w_q = f_q \quad \text{on } (-\infty, a),$$

subject to boundedness as $x \to -\infty$ and the Robin condition

$$\mathcal{B}_a w_q = -\lambda.$$

By Lemma 5.1, $\psi_q$ and $\chi_q$ are two linearly independent solutions of the homogeneous equation $L_q u = 0$, where $\psi_q$ is the unique homogeneous solution that remains bounded as $x \to -\infty$, and $\chi_q$ satisfies $\mathcal{B}_a \chi_q = 0$. These two solutions are therefore adapted to the two boundary requirements of the problem.

Let

$$W_q(x) := \psi_q(x)\chi_q'(x) - \psi_q'(x)\chi_q(x)$$

denote their Wronskian. Since $\psi_q$ and $\chi_q$ are linearly independent, one has $W_q(x) \neq 0$ on $(-\infty, a)$. We first solve the associated homogeneous boundary-value problem

$$L_q w = 0, \quad w \text{ bounded as } x \to -\infty, \quad \mathcal{B}_a w = -\lambda.$$

Because $\psi_q$ is the unique bounded homogeneous solution on the left, such a solution must be of the form

$$w(x) = A_q \psi_q(x)$$

for some constant $A_q$. Imposing the Robin condition gives

$$A_q \mathcal{B}_a \psi_q = -\lambda.$$

Hence

$$A_q = -\frac{\lambda}{\mathcal{B}_a \psi_q}.$$

By Lemma 5.1, $\mathcal{B}_a \psi_q = \frac{\sigma\sqrt{b}}{\sqrt{2}} e^{p(a)} D_{v_q+1}(z(a))$, and therefore

$$A_q = -\frac{\sqrt{2}\,\lambda}{\sigma\sqrt{b}} \frac{e^{-p(a)}}{D_{v_q+1}(z(a))}.$$

It follows that

$$w_q^{(0)}(x) = A_q \psi_q(x) = -\frac{\sqrt{2}\,\lambda}{\sigma\sqrt{b}} \exp(p(x) - p(a)) \frac{D_{v_q}(z(x))}{D_{v_q+1}(z(a))},$$

which is exactly the expression in (60).

We now incorporate the inhomogeneous term $f_q$. Standard Green-function theory for second-order linear equations shows that the unique solution of $L_q w = f_q$ which is bounded as $x \to -\infty$ and satisfies the homogeneous Robin condition $\mathcal{B}_a w = 0$ is given by

$$w(x) = \int_{-\infty}^{a} \Gamma_q(x,y) f_q(y)\, dy,$$

where the Green kernel is

$$\Gamma_q(x,y) = \frac{2}{\sigma^2 W_q(y)} \begin{cases} \psi_q(x) \chi_q(y), & x \leq y < a, \\ \psi_q(y) \chi_q(x), & y \leq x < a. \end{cases}$$

Indeed, for fixed $y$, the function $x \mapsto \Gamma_q(x,y)$ solves the homogeneous equation $L_q u = 0$ away from $x = y$, is bounded on the left because it is proportional to $\psi_q(x)$ for $x \leq y$, satisfies the homogeneous Robin condition at $a$ because it is proportional to $\chi_q(x)$ for $x \geq y$, and has the correct derivative jump at $x = y$ determined by the factor $2/\left(\sigma^2 W_q(y)\right)$.

Combining the homogeneous solution $w_q^{(0)}$ with this Green representation, we conclude that the full solution of

$$L_q w_q = f_q, \qquad \mathcal{B}_a w_q = -\lambda, \qquad w_q \text{ bounded as } x \to -\infty,$$

is

$$w_q(x) = w_q^{(0)}(x) + \int_{-\infty}^{a} \Gamma_q(x,y) f_q(y)\, dy.$$

Finally, substitute the specific forcing term

$$f_q(y) = \eta q \int_y^a w_q(u)\, du.$$

This gives

$$w_q(x) = w_q^{(0)}(x) + \eta q \int_{-\infty}^{a} \Gamma_q(x,y) \left( \int_y^a w_q(u)\, du \right) dy, \qquad x < a,$$

which is exactly the Green–Volterra representation in (59). □

Theorem 5.2 shows that the discounted problem remains explicitly analyzable even though it no longer collapses to a closed parabolic-cylinder formula: the nonlocal effect of discounting is isolated entirely in a Volterra perturbation of the Weber operator.

*Remark 5.3* When $q = 0$, the Volterra term vanishes identically. The Green–Volterra representation therefore collapses to the explicit formula $w_0(x) = w_0^{(0)}(x)$, which is exactly the starting point for the parabolic-cylinder expression obtained below for the undiscounted overshoot probability.

**Theorem 5.4** (Small-$q$ expansion of the discounted overshoot transform).

*Assume $\beta < 0$ and fix $x < a$. Suppose that $\mathbb{E}_x[\tau_a^2] < \infty$. Then, as $q \downarrow 0$,*

$$G_q(x) = G_0(x) - qM(x) + R_q(x), \qquad M(x) := \mathbb{E}_x[\tau_a \mathbf{1}_{J_+}], \tag{62}$$

*with remainder bounded by*

$$|R_q(x)| \le \frac{q^2}{2} \mathbb{E}_x[\tau_a^2 \mathbf{1}_{J_+}]. \tag{63}$$

*Equivalently,*

$$\frac{G_0(x) - G_q(x)}{q} \to M(x) \quad \text{as } q \downarrow 0. \tag{64}$$

*Moreover, if*

$$m_q(x) := \frac{G_0(x) - G_q(x)}{q}, \tag{65}$$

*then each $m_q$ is the unique bounded classical solution of*

$$\frac{1}{2}\sigma^2 m_q''(x) + (\alpha - bx)m_q'(x) + \lambda \int_0^{a-x}\bigl(m_q(x+z) - m_q(x)\bigr)\eta e^{-\eta z}\,dz \\ = -G_q(x), \qquad x < a, \tag{66}$$

*with boundary condition $m_q(a^-) = 0$.*

*Proof.* Fix $x < a$, and recall that $G_q(x) = \mathbb{E}_x[e^{-q\tau_a}\mathbf{1}_{J_+}]$, $G_0(x) = \mathbb{E}_x[\mathbf{1}_{J_+}] = \mathbb{P}_x(J_+)$.

Hence

$$G_0(x) - G_q(x) = \mathbb{E}_x[(1 - e^{-q\tau_a})\mathbf{1}_{J_+}].$$

We first derive the small-$q$ expansion. For every $u \ge 0$,

$$0 \le 1 - e^{-u} \le u,$$

and the Taylor remainder formula yields

$$|1 - e^{-u} - u| \leq \frac{u^2}{2}.$$

Applying these inequalities with $u = q\tau_a$, we obtain

$$0 \leq \frac{1 - e^{-q\tau_a}}{q} \leq \tau_a,$$

and therefore

$$0 \leq \frac{G_0(x) - G_q(x)}{q} \leq \mathbb{E}_x[\tau_a \mathbf{1}_{J_+}].$$

Moreover,

$$|1 - e^{-q\tau_a} - q\tau_a| \leq \frac{q^2 \tau_a^2}{2}.$$

Multiplying by $\mathbf{1}_{J_+}$ and taking expectations gives

$$|G_q(x) - G_0(x) + q\, \mathbb{E}_x[\tau_a \mathbf{1}_{J_+}]| \leq \frac{q^2}{2} \mathbb{E}_x[\tau_a^2 \mathbf{1}_{J_+}].$$

Thus, with $M(x) := \mathbb{E}_x[\tau_a \mathbf{1}_{J_+}]$, we have

$$G_q(x) = G_0(x) - qM(x) + R_q(x),$$

where the remainder satisfies $|R_q(x)| \leq \frac{q^2}{2} \mathbb{E}_x[\tau_a^2 \mathbf{1}_{J_+}]$. Dividing by $q$ and letting $q \downarrow 0$ yields (64), which proves the first part of the theorem.

Now define

$$m_q(x) := \frac{G_0(x) - G_q(x)}{q}.$$

We derive the boundary-value problem satisfied by $m_q$. For $q > 0$, the discounted transform $G_q$ satisfies

$$\frac{1}{2}\sigma^2 G_q''(x) + (\alpha - bx)G_q'(x) + \lambda \int_0^{a-x} \big(G_q(x+z) - G_q(x)\big) \eta e^{-\eta z}\, dz - qG_q(x) + \lambda e^{-\eta(a-x)}$$
$$= 0.$$

At $q = 0$, the same equation reads

$$\frac{1}{2}\sigma^2 G_0''(x) + (\alpha - bx)G_0'(x) + \lambda \int_0^{a-x} \big(G_0(x+z) - G_0(x)\big)\eta e^{-\eta}\, dz + \lambda e^{-\eta(a-x)} = 0.$$

Subtracting the second equation from the first, the source term $\lambda e^{-\eta(a-x)}$ cancels, and we obtain

$$\frac{1}{2}\sigma^2(G_q'' - G_0'')(x) + (\alpha - bx)(G_q' - G_0')(x)$$
$$+ \lambda \int_0^{a-x} \bigl((G_q - G_0)(x+z) - (G_q - G_0)(x)\bigr)\eta e^{-\eta z}\, dz - q\, G_q(x) = 0.$$

Dividing by $q$ and using $\frac{G_q(x) - G_0(x)}{q} = -m_q(x)$ gives (66).

The boundary condition follows immediately from the fact that
$$G_q(a^-) = 0 \quad \text{and} \quad G_0(a^-) = 0,$$
namely
$$m_q(a^-) = \frac{G_0(a^-) - G_q(a^-)}{q} = 0.$$

It remains to prove uniqueness in the class of bounded classical solutions. Let $m_q^{(1)}$ and $m_q^{(2)}$ be two such solutions, and set $v := m_q^{(1)} - m_q^{(2)}$. Then $v$ is bounded, satisfies $v(a^-) = 0$, and solves the homogeneous equation
$$\frac{1}{2}\sigma^2 v''(x) + (\alpha - bx)v'(x) + \lambda \int_0^{a-x} \bigl(v(x+z) - v(x)\bigr)\eta e^{-\eta z}\, dz = 0, \quad x < a.$$

Applying Dynkin's formula to $v(X_{t\wedge \tau_a})$, we obtain
$$v(x) = \mathbb{E}_x\bigl[v(X_{t\wedge \tau_a})\bigr].$$

Now let $t \to \infty$. Since $v$ is bounded and, by Theorem 6.1, $\mathbb{P}_x(\tau_a < \infty) = 1$ in the mean-reverting case, the terminal variable converges to $v(a^-) = 0$. Hence
$$v(x) = 0 \quad \text{for all } x < a.$$

Therefore, $m_q^{(1)} \equiv m_q^{(2)}$, which proves uniqueness. □

*Remark 5.5* The coefficient $M(x) = \mathbb{E}_x[\tau_a \mathbf{1}_{J_+}]$ is the first discounted correction to the overshoot probability. Thus the small-$q$ regime couples the overshoot mode not only to its probability $G_0(x)$ but also to the overshoot-conditioned time scale through the first moment of $\tau_a$. Differentiating the identity $G_q(x) = -\int_x^a w_q(y)\, dy$ suggests that the first derivative of $w_q$ at $q = 0$, when it exists, should solve the inhomogeneous Weber-type equation obtained by differentiating (50); the preceding theorem identifies the corresponding source term at the level of $G_q$ itself.

## Section 6 — Undiscounted overshoot probability and creeping

This section focuses on the undiscounted overshoot probability $G_0(x) = \mathbb{P}_x(J_+)$.

**Theorem 6.1** (Explicit overshoot mode probability in the mean-reverting case).

*Assume $\beta < 0$ and define*

$$b := -\beta > 0, \qquad v := -1 - \frac{\lambda}{b}, \qquad z(x) := \frac{\sqrt{2}}{\sigma\sqrt{b}}\left(-bx + \alpha + \frac{1}{2}\eta\sigma^2\right),$$

$$p(x) := -\frac{\alpha x}{\sigma^2} + \frac{bx^2}{2\sigma^2} + \frac{\eta x}{2}.$$

*Then*

$$G_0'(x) = -\frac{\sqrt{2}\lambda}{\sigma\sqrt{b}} \exp(p(x) - p(a)) \frac{D_v(z(x))}{D_{v+1}(z(a))}, \qquad x < a, \tag{67}$$

*and therefore*

$$G_0(x) = \frac{\sqrt{2}\lambda}{\sigma\sqrt{b}\, D_{v+1}(z(a))} \int_x^a \exp(p(y) - p(a)) D_v(z(y))\, dy. \tag{68}$$

*Moreover,*

$$\mathbb{P}_x(\tau_a < \infty) = 1, \qquad x < a, \tag{69}$$

*and hence*

$$\mathbb{P}_x(C_0) = 1 - G_0(x). \tag{70}$$

*Proof.*

Assume $\beta < 0$, and write $b := -\beta > 0$. At $q = 0$, the reduction obtained earlier shows that $w(x) := G_0'(x)$ satisfies the second-order ODE

$$\frac{1}{2}\sigma^2 w''(x) + \left(\alpha - bx - \frac{1}{2}\eta\sigma^2\right)w'(x) + (-b - \eta(\alpha - bx) - \lambda)w(x) = 0, \qquad x < a.$$

Introduce the gauge transform

$$w(x) = e^{p(x)}\phi(x), \qquad p(x) := -\frac{\alpha x}{\sigma^2} + \frac{bx^2}{2\sigma^2} + \frac{\eta x}{2},$$

and the new variable

$$z(x) := \frac{\sqrt{2}}{\sigma\sqrt{b}}\left(-bx + \alpha + \frac{1}{2}\eta\sigma^2\right).$$

A direct calculation shows that $\phi$ then satisfies

$$\phi_{zz}(z) + \left(\nu + \frac{1}{2} - \frac{z^2}{4}\right)\phi(z) = 0, \qquad \nu := -1 - \frac{\lambda}{b}.$$

This is the parabolic-cylinder equation. Its solutions are linear combinations of $D_\nu(z)$ and a second independent branch that grows exponentially as $z \to +\infty$.

Now $z(x) \to +\infty$ as $x \to -\infty$. Since $0 \le G_0(x) \le 1$, the overshoot probability is bounded on $(-\infty, a)$. This excludes the exponentially growing branch, and therefore

$$G_0'(x) = C\, e^{p(x)} D_\nu(z(x))$$

for some constant $C$. To determine $C$, use the compatibility condition obtained earlier from the OIDE/ODE equivalence:

$$\frac{1}{2}\sigma^2 G_0''(a^-) + (\alpha - ba)G_0'(a^-) + \lambda = 0.$$

Differentiate the expression for $G_0'$. Since

$$z'(x) = -\frac{\sqrt{2b}}{\sigma}, \quad D_\nu'(z) = \frac{z}{2}D_\nu(z) - D_{\nu+1}(z), \quad \frac{1}{2}\sigma^2 p'(a) + (\alpha - ba) = \frac{\sigma\sqrt{b}}{2\sqrt{2}}z(a),$$

substituting into the compatibility condition shows that the terms involving $D_\nu(z(a))$ cancel, leaving

$$C\, e^{p(a)} \frac{\sigma\sqrt{b}}{\sqrt{2}} D_{\nu+1}(z(a)) + \lambda = 0.$$

Hence

$$C = -\frac{\sqrt{2}\,\lambda}{\sigma\sqrt{b}} \frac{e^{-p(a)}}{D_{\nu+1}(z(a))},$$

and (67) is obtained. Since $G_0(a^-) = 0$, integrating (67) from $x$ to $a$ yields

$$G_0(x) = -\int_x^a G_0'(y)\, dy,$$

and (68) ensues. It remains to prove that

$$\mathbb{P}_x(\tau_a < \infty) = 1.$$

Let $Y$ be the Ornstein–Uhlenbeck diffusion driven by the same Brownian motion,

$$dY_t = (\alpha - bY_t)\, dt + \sigma\, dW_t, \qquad Y_0 = x.$$

The variation-of-constants formula gives

$$X_t = Y_t + \int_0^t e^{-b(t-s)}\, dJ_s.$$

Since $J$ is increasing, the integral term is nonnegative, so

$$X_t \geq Y_t \quad \text{for all } t \geq 0.$$

Next, the Ornstein–Uhlenbeck process has the explicit representation

$$Y_t = \frac{\alpha}{b} + e^{-bt}\left(x - \frac{\alpha}{b}\right) + \sigma e^{-b} \tilde{B}_{\rho(t)}, \qquad \rho(t) := \frac{e^{2b} - 1}{2b}.$$

for a Brownian motion $\tilde{B}$. Since

$$e^{-b} = \frac{1}{\sqrt{1 + 2b\rho(t)}},$$

we may rewrite the fluctuation term as

$$\sigma e^{-bt} \tilde{B}_{\rho(t)} = \frac{\sigma \tilde{B}_s}{\sqrt{1 + 2bs}}, \qquad s = \rho(t).$$

By the law of the iterated logarithm,

$$\limsup_{s \to \infty} \frac{\tilde{B}_s}{\sqrt{2s \log\log s}} = 1 \quad \text{a.s.}$$

Hence, along a subsequence,

$$\tilde{B}_s \sim \sqrt{2s \log\log s},$$

and therefore

$$\frac{\tilde{B}_s}{\sqrt{1 + 2bs}} \sim \sqrt{\frac{\log\log s}{b}} \to +\infty.$$

It follows that

$$\limsup_{t \to \infty} \sigma e^{-bt} \tilde{B}_{\rho(t)} = +\infty \quad \text{a.s.}$$

and thus

$$\sup_{t \geq 0} Y_t = +\infty \quad \text{a.s.}$$

Consequently,

$$\tau_a^Y := \inf\{t \geq 0 : Y_t \geq a\} < \infty \quad \text{a.s.}$$

Since $X_t \geq Y_t$ for all $t$, one has $\tau_a \leq \tau_a^Y$ and therefore (69) holds.

Finally, for this model the earlier dichotomy result yields

$$\{\tau_a < \infty\} = C_0 \,\dot\cup\, J_+ \quad \text{a.s.}$$

Since $\mathbb{P}_x(\tau_a < \infty) = 1$ and $G_0(x) = \mathbb{P}_x(J_+)$, (70) ensues. □

*Remark 6.2* (Boundary asymptotics).

The following formula gives the one-sided boundary expansion:

$$G_0(x) = \frac{\sqrt{2}\,\lambda}{\sigma\sqrt{b}} \frac{D_\nu(z(a))}{D_{\nu+1}(z(a))} (a-x) + o(a-x), \quad x \uparrow a. \tag{71}$$

Consequently,

$$\mathbb{P}_x(C_0) = 1 - \frac{\sqrt{2}\,\lambda}{\sigma\sqrt{b}} \frac{D_\nu(z(a))}{D_{\nu+1}(z(a))} (a-x) + o(a-x), \quad x \uparrow a. \tag{72}$$

Thus continuous contact dominates arbitrarily close to the barrier, while the jump-overshoot mode enters through an explicit slope depending on $(\alpha, \beta, \sigma, \lambda, \eta, a)$.

*Remark 6.3* When $\beta = 0$, the affine equation reduces to the constant-drift benchmark, in which the third-order ODE has constant coefficients. The mean-reverting case treated above is analytically different: the state dependence in the drift removes spatial homogeneity but still leaves the overshoot mode explicitly tractable through the parabolic-cylinder reduction.

# 7 Conclusion

The natural starting point for first passage through a continuous barrier by a càdlàg process is the fourfold decomposition. This decomposition is purely pathwise, but it has a clear random-time interpretation: the left-contact component always defines an accessible stopping time and, under Assumption 2.5, the canonical running-supremum sequence upgrades it to a predictable one, while under a no-predictable-jumps condition the gap-from-the-left component defines a totally inaccessible stopping time. The familiar contact-versus-overshoot dichotomy appears only after excluding jumps away from the graph and exact hit by jump.

In the semimartingale setting, Assumption 2.5 admits a sharp compensator criterion, and the jump-driven crossing modes admit exact characteristic formulas. The affine analysis of Section 6 shows that these formulas can be combined with Dynkin's formula to produce a genuine boundary-value problem: for upward exponential jumps, the overshoot mode satisfies an integro-differential equation whose differentiated third-order ODE is equivalent to the original problem only after one retains the boundary compatibility condition. For the discounted problem, the affine OIDE also admits a verification and uniqueness theorem in the class of bounded classical solutions, and in the mean-reverting case the discounted overshoot transform is further reduced to an explicit Green–Volterra representation built from parabolic-cylinder functions together with a first-order small-$q$ expansion whose leading correction is the overshoot-time moment. The undiscounted problem then collapses to a closed parabolic-cylinder formula for the overshoot probability, and the creeping probability follows by complementarity.

Several directions remain open. One may sharpen the discounted analysis of $F_q$, $G_q$, and $H_q$, treat nonconstant continuous barriers in affine models, or replace the exponential jump law by phase-type or hyperexponential distributions, where one expects a higher-order local system after augmenting the state space. More broadly, the framework provides a bridge between fluctuation theory, the general theory of stochastic processes, and the structural intuition of visible versus surprise occurrence.

# References

[1] A. Aksamit and M. Jeanblanc. Enlargement of Filtration with Finance in View. Springer, 2017.

[2] L. Alili and A. E. Kyprianou. Some remarks on first passage of Lévy processes, the American put and pasting principles. *Ann. Appl. Probab.* **15** (2005), 2062–2080.

[3] R. F. Bass. The measurability of hitting times. *Electron. Commun. Probab.* **15** (2010), 99–105.

[4] F. Black and J. C. Cox. Valuing corporate securities: Some effects of bond indenture provisions. J. Finance 31 (1976), 351–367.

[5] N. Cai. On first passage times of a hyper-exponential jump diffusion process. Oper. Res. Lett. 37 (2009), 127–134.

[6] U. Çetin, R. Jarrow, P. Protter, and Y. Yildirim. Modeling credit risk with partial information. *Ann. Appl. Probab.* **14** (2004), 1167–1178.

[7] L. Chaumont and T. Pellas. Creeping of Lévy processes through curves. *Ann. Appl. Probab.* **33** (2023), 2609–2642.

[8] Y.-T. Chen, Y.-C. Sheu, and M.-C. Chang. A note on first passage functionals for hyper-exponential jump-diffusion processes. Electron. Commun. Probab. 18 (2013), 1–8.

[9] R. A. Doney and A. E. Kyprianou. Overshoots and undershoots of Lévy processes. *Ann. Appl. Probab.* **16** (2006), 91–106.

[10] R. J. Elliott, M. Jeanblanc, and M. Yor. On models of default risk. *Math. Finance* **10** (2000), 179–195.

[11] K. Giesecke. Default and information. J. Econom. Dynam. Control 30 (2006), 2281–2303.

[12] J. He, S. Wang, and J. Yan. *Semimartingale Theory and Stochastic Calculus*. CRC Press, 1992.

[13] M. JACOBSEN AND A. T. JENSEN. Exit times for a class of piecewise exponential Markov processes with two-sided jumps. *Stochastic Process. Appl.* **117** (2007), 1330–1356.

[14] M. JEANBLANC AND S. VALCHEV. Partial information, default hazard process, and default-risky bonds. *Int. J. Theor. Appl. Finance* **8** (2005), 807–838.

[15] S. G. KOU. A jump-diffusion model for option pricing. *Management Sci.* **48** (2002), 1086–1101.

[16] S. G. KOU AND H. WANG. First passage times of a jump diffusion process. *Adv. in Appl. Probab.* **35** (2003), 504–531.

[17] S. G. KOU AND H. WANG. Option pricing under a double exponential jump diffusion model. *Management Sci.* **50** (2004), 1178–1192.

[18] A. E. KYPRIANOU. *Fluctuations of Lévy Processes with Applications*, 2nd ed. Springer, 2014.

[19] A. E. KYPRIANOU, J. C. PARDO, AND V. RIVERO. Exact and asymptotic $n$-tuple laws at first and last passage. *Ann. Appl. Probab.* **20** (2010), 522–564.

[20] R. C. MERTON. Option pricing when underlying stock returns are discontinuous. *J. Financial Econom.* **3** (1976), 125–144.

[21] A. NIKEGHBALI. An essay on the general theory of stochastic processes. *Probab. Surv.* **3** (2006), 345–412.

[22] A. A. NOVIKOV. Martingales and first-passage times for Ornstein–Uhlenbeck processes with a jump component. *Theory Probab. Appl.* **48** (2004), 288–303.

[23] M. SONG. On the ruin problem in the renewal risk processes perturbed by diffusion. *Ann. Inst. Stat. Math.* **61** (2009), 135–153.

[24] A. W. VAN DER VAART. *Martingales, Diffusions and Financial Mathematics*. Springer, 2014.

[25] C. YIN, Y. WEN, Z. ZONG, AND Y. SHEN. The first passage time problem for mixed-exponential jump processes with applications in insurance and finance. Abstr. Appl. Anal. 2014 (2014), 571724.

[26] C. ZHOU. The term structure of credit spreads with jump risk. J. Banking Finance 25 (2001), 2015–2040.